\documentclass[11pt]{amsart}
\usepackage{amssymb,amscd,latexsym,epsfig,xypic,hyperref}
\usepackage[mathscr]{euscript}

\DeclareRobustCommand{\SkipTocEntry}[4]{}
\makeatletter
\newcommand\@dotsep{4.5}
\def\@tocline#1#2#3#4#5#6#7{\relax
  \ifnum #1>\c@tocdepth 
  \else
    \par \addpenalty\@secpenalty\addvspace{#2}%
    \begingroup \hyphenpenalty\@M
    \@ifempty{#4}{%
      \@tempdima\csname r@tocindent\number#1\endcsname\relax
    }{%
      \@tempdima#4\relax
    }%
    \parindent\z@ \leftskip#3\relax \advance\leftskip\@tempdima\relax
    \rightskip\@pnumwidth plus1em \parfillskip-\@pnumwidth
    #5\leavevmode\hskip-\@tempdima #6\relax
    \leaders\hbox{$\m@th
      \mkern \@dotsep mu\hbox{.}\mkern \@dotsep mu$}\hfill
    \hbox to\@pnumwidth{\@tocpagenum{#7}}\par
    \nobreak
    \endgroup
  \fi}
\makeatother 

\DeclareFontFamily{OT1}{rsfs}{}
\DeclareFontShape{OT1}{rsfs}{n}{it}{<-> rsfs10}{}
\DeclareMathAlphabet{\curly}{OT1}{rsfs}{n}{it}

\DeclareFontFamily{OT1}{arrow}{\hyphenchar\font45 }
\DeclareFontShape{OT1}{arrow}{m}{n}{<7><8><10><12><13.82><16.59><19.907><23.89><28.66><34.4><41.28>cmr8}{}
\DeclareSymbolFont{arrowsym}{OT1}{arrow}{m}{n}
\DeclareMathSymbol{\arroweq}{\mathrel}{arrowsym}{"3D}

\newcommand\PP{\mathbb P}
\newcommand\C{\mathbb C}
\newcommand\Q{\mathbb Q}

\newcommand\Z{\mathbb Z}

\renewcommand\O{\mathcal O}

\newcommand\Mb{\,\overline{\!M}}

\newcommand{\Rt}[1]{\stackrel{#1\,}{\longrightarrow}}
\newcommand{\RT}[2]{\xymatrix@C=#1pt{\ar[r]^{#2}&}}
\newcommand\To{\longrightarrow}
\newcommand\into{\hookrightarrow}
\newcommand\INTO{\ \ar@{^(->}[r]<-.2ex>}

\renewcommand\_{^{}_}

\newfont{\bigtimesfont}{cmsy10 scaled \magstep5}
\newcommand{\bigtimes}{\mathop{\lower0.9ex\hbox{\bigtimesfont\symbol2}}}
\newcommand\udot{^\bullet}

\newcommand\vir{\operatorname{vir}}
\newcommand\vd{\operatorname{vd}}
\newcommand\red{\operatorname{red}}

\newcommand\Hom{\operatorname{Hom}}

\newcommand\Ext{\operatorname{Ext}}

\newcommand\Aut{\operatorname{Aut}}

\newcommand\Spec{\operatorname{Spec}\,}
\newcommand\Hilb{\operatorname{Hilb}}
\newcommand\Bl{\operatorname{Bl}}

\renewcommand\;{\hspace{.6pt}}

\newcommand\beq[1]{\begin{equation}\label{#1}}
\newcommand\eeq{\end{equation}}
\newcommand\beqa{\begin{eqnarray*}}
\newcommand\eeqa{\end{eqnarray*}}

\makeatletter \@addtoreset{equation}{section} \makeatother

\newtheorem{thm}[equation]{Theorem}

\newtheorem{prop}[equation]{Proposition}

\title{Notes on the proof of the KKV conjecture}
\author[R. Pandharipande and R. P. Thomas]{Rahul Pandharipande and Richard P. Thomas}

\begin{document}
\begin{abstract} \noindent
The Katz-Klemm-Vafa conjecture 
expresses the Gromov-Witten theory of K3 surfaces 
(and K3-fibred 3-folds in fibre classes) 
in terms of modular forms. Its recent proof gives the first non-toric geometry in dimension greater than 1 where Gromov-Witten theory is exactly solved in all genera.

We survey the various steps in the proof. 
The MNOP correspondence
and a new Pairs/Noether-Lefschetz correspondence for K3-fibred 3-folds 
transform the Gromov-Witten problem into a calculation of the full stable pairs theory 
of a local K3-fibred 3-fold. The stable pairs
calculation is then carried out via degeneration, 
localisation, vanishing results, and new multiple cover formulae.
\end{abstract}

\maketitle


\setcounter{tocdepth}{2} 
\tableofcontents

\setcounter{section}{-1}
\section{Introduction}

Gromov-Witten theory provides a way to count curves in algebraic varieties and symplectic manifolds. The theory is essentially nonlinear, and computations 
are very difficult. Furthermore, degenerate contributions make the link to 
the enumeration of curves rather
opaque. Complete calculations have been possible only for targets which are either toric or have dimension at most $1$.

Nevertheless, for K3 surfaces and curve classes in the fibres of K3-fibred 
3-folds, the full Gromov-Witten theory was conjecturally described in \cite{KKV, MP} in terms of modular forms. This paper is a survey of the conjecture and 
our recent proof \cite{PTKKV}.

K3 surfaces provide a basic
example of the practical use of the MNOP correspondence 
 for computations. 
The recent proof of the  correspondence
for many Calabi-Yau 3-folds \cite{PaPix} converts a calculation in Gromov-Witten theory 
into a slightly more linear calculation in sheaf theory. While the sheaf theory questions are still hard, 
we will see, for K3 surfaces, that they are solvable. At present, no direct approach to the integrals in
Gromov-Witten theory is available.
\medskip

In brief, Sections \ref{GW} to \ref{MNOPsec} review the relevant theory and set up the problem. Sections \ref{localMNOP} to \ref{invert} prove a local MNOP conjecture by combining the global MNOP correspondence
 of \cite{PaPix} with a new Pairs/Noether-Lefschetz correspondence. Finally Sections \ref{EC} to \ref{mult} compute the resulting sheaf theory problem in full. An outline of the proof goes as follows, numbered by Section.

\begin{enumerate}
\setcounter{enumi}{0}
\item A brief review of stable maps, Gromov-Witten invariants, multiple covers and degenerate contributions. For a fuller introduction to various curve counting theories and their
interrelationships -- especially the MNOP conjecture -- we refer the reader to \cite{13.5}.
\item The Gromov-Witten invariants of 3-folds are recast in terms of BPS numbers, which are conjecturally integers.
\item We describe the relevant Gromov-Witten invariants of K3 surfaces from a 2-dimensional point of view -- via the {\em reduced virtual class} on the K3 surface $S$ -- and from a 3-dimensional point of view via a local (algebraic approximation to the) {\em twistor 3-fold} $T$.
\item We state the KKV conjecture and discuss its predictions for multiple covers.
\item We give a brief review of the sheaf theory we use to count curves, called {\em stable pairs}.
\item A brief review of the MNOP conjecture relating Gromov-Witten and stable pair invariants for 3-folds.
\item We discuss how to derive an MNOP correspondence for the local twistor 
3-fold $T\to\Delta$ from the proof of the MNOP conjecture for projective Calabi-Yau 3-folds $X$ in \cite{PaPix}. The proof will occupy Sections \ref{JLi} to \ref{invert}.
\item We describe relative geometries and the degeneration formula for
counting invariants.
\item This is applied to the deformation to the normal cone of $S$ inside a K3-fibration.
\item The result is a surprising description of the \emph{connected} version of the stable pair invariants.
\item This allows us to prove a 
Pairs/Noether-Lefschetz correspondence for K3-fibred 3-folds. The stable pairs theory of a K3-fibred 3-fold is described in terms of the stable pairs theory of K3 surfaces and Noether-Lefschetz numbers counting the number of K3 fibres for which our curve class is algebraic.
\item We choose a convenient projective K3-fibred 3-fold $X$ for which
 the MNOP correspondence
 is known \cite{PaPix}, and the 
Pairs/Noether-Lefschetz correspondence on $X$ is \emph{invertible}: the stable pair invariants of a K3 surface can be recovered from the collection of Noether-Lefschetz numbers and stable pair invariants of $X$. Together with the Gromov-Witten/Noether-Lefschetz correspondence of \cite{MP} we deduce the required local MNOP conjecture for K3 surfaces.
\item The entire problem is now translated into one of computing the stable pairs theory of the local 3-fold $S\times\C$. The multiple covers of Gromov-Witten theory have become stable pairs on $S\times\{0\}$ which are scheme-theoretically thickened in the $\C$-direction.

\item We describe a critical advantage of stable pairs theory -- its symmetric obstruction theory.
This defines a linear functional on the obstruction theory which forces many of the invariants to vanish. This vanishing result simplifies the multiple cover structure of stable pairs theory considerably.
\item Finally we compute these multiple covers. Calculations in \cite{KT1,KT2} show they have the remarkable topological properties required for the KKV conjecture. 
This reduces the problem to the primitive case where previous calculations of Kawai-Yoshioka \cite{KY,MPT} give the KKV formula.
\end{enumerate} \medskip

For the reader wishing to use this survey as a guide to the full proofs in \cite{PTKKV}, we give the following approximate correspondence between the various sections. \medskip

\begin{center}
\begin{tabular}{|c||c|c|c|c|c|c|c|c|c|}
        \hline
\textbf{}
Section of the paper \cite{PTKKV} &2&3&4&5&6&7&8 \\
\hline Relevant sections of this paper
&12&7&15&13, 14&9, 10&15&11 \\ \hline
\end{tabular}
\end{center}
 

\section{Gromov-Witten theory} \label{GW}
We refer to \cite{13.5} for an introduction to curve counting on nonsingular
projective varieties
$X$ in class $\beta\in H_2(X,\Z)$.
We give here a brief review of only what is necessary to state the KKV conjecture. 

A \emph{stable map} is the data of
\begin{itemize}
\item an algebraic map $f\colon C\to X$, where
\item $C$ is a connected algebraic curve, at worst nodal, and 
\item $|\!\Aut(f)|<\infty$.
\end{itemize} 
Here, $\Aut(f)$ is the group of automorphisms of $C$ which fix $f$. Representing the associated moduli problem is a moduli space $$\Mb_g(X,\beta)$$ of stable maps whose domain curve $C$ has arithmetic genus $g$. It is a
\emph{compact} Deligne-Mumford stack with finite stabiliser groups due to the third condition above.

Although $\Mb_g(X,\beta)$ is usually singular -- so the dimension of the deformation space\footnote{When $f$ is an embedding the complex becomes quasi-isomorphic to the conormal bundle $N^*_C$ to $C\into X$, so the deformation space becomes the more familiar $H^0(N_C)$.} $\Hom\!\big(\big\{f^*\Omega_X\to\Omega_C\big\},\O_C\big)$ of a stable map $f$ is unpredictable -- there is a natural obstruction theory $\Ext^1\!\big(\big\{f^*\Omega_X\to\Omega_C\big\},\O_C\big)$ such that the difference in dimensions
$$
\dim\Hom\!\big(\big\{f^*\Omega_X\to\Omega_C\big\},\O_C\big)-\dim\Ext^1\!
\big(\big\{f^*\Omega_X\to\Omega_C\big\},\O_C\big)
$$
is the topological constant
\beq{vd}
\vd=\int_\beta\!c\_1(X)\ +\ (1-g)(\dim X-3),
\eeq
called the virtual dimension. The components of the
moduli space always have dimension at least $\vd$. 
There exists a  virtual moduli cycle
\begin{equation}\label{jj33}
\big[\Mb_g(X,\beta)\big]^{\vir}\in A_{\vd}\big(\Mb_g(X,\beta)\big)\To H_{2\!\vd}\big(\Mb_g(X,\beta)\big)
\end{equation}
which is the usual fundamental cycle of $\Mb_g(X,\beta)$ when the dimension of the moduli space equals $\vd$ \eqref{vd}.

Integrating cohomology classes over the  cycle \eqref{jj33}
gives Gromov-Witten invariants. Our main interest will be in Calabi-Yau 3-folds, where the formula \eqref{vd} gives $\vd=0$. Then, the only Gromov-Witten invariant is the degree of the virtual cycle,
\beq{CYGW}
N_{g,\beta}(X)=\int_{[\Mb_g(X,\beta)]^{\vir}}1\ \in\,\Q\, ,
\eeq
which provides a virtual count of the {\em number of curves of degree $\beta$ and genus $g$ in $X$}.
The result
is a rational number because of the finite automorphisms of stable maps.
The integral $N_{g,\beta}(X)$ 
 is invariant under deformations of $X$.

\subsection*{Multiple covers and degenerate contributions}
It is notoriously hard to calculate and interpret \eqref{CYGW} due to multiple covers and 
{\em degenerate contributions}.

The simplest examples of these two phenomena come from considering a nonsingular isolated $(-1,-1)$ rational curve
$$
C\cong\PP^1\subset X
$$
with fundamental class $\beta$. For simplicity, assume there are no other curves in class $d\beta$ for any $d$ (other than the obvious covers and thickenings of $C$ whose set-theoretic image is just $C$ itself).

The embedding of $C$ is a stable map which gives
\beq{1}
N_{0,\beta}(X)=1,
\eeq
as expected. However, degree $d$ covers of $C$ also contribute in class $d\beta$. Since they have finite automorphisms -- the deck transformations of the cover -- the contribution is rational. In genus 0, for instance, 
the answer is given by the Aspinwall-Morrison formula
\beq{AspMor}
N_{0,d\beta}\ =\ \frac1{d^3}\,.
\eeq
Things are even worse in higher genus $g\ge1$. Glue any genus $g$ curve $\Sigma$ to $C$ at a point $p$ and map the resulting nodal curve $\Sigma\cup_p\PP^1$ to $C\subset X$ by contracting $\Sigma$. The resulting stable maps contribute to $N_{g,\beta}(X)$ through integrals over the moduli space $\Mb_{g,1}$ of abstract pointed stable curves $(\Sigma,p)$. These \emph{degenerate contributions} to the Gromov-Witten invariant \eqref{CYGW} are computed in \cite{FP,Pdegen}; in more general situations they contribute significant complication to Gromov-Witten theory.

\section{BPS reformulation for 3-folds} \label{GVsec}

In the above example, there is an obvious integer underlying the rational numbers $N_{g,d\beta}(X)$: it is $N_{0,\beta}(X)=1$ \eqref{1}, from which all the others can be derived by universal formulae.

Gopakumar and Vafa predicted that this phenomenon should hold in general for Calabi-Yau 3-folds.\footnote{An extension to all 3-folds is given in \cite{PICM}.} The underlying integer counts are called the BPS invariants $n_{g,\beta}(X)$. They proposed a sheaf-theoretic definition which has yet to be made mathematically precise in general.\footnote{But see \cite[Appendix A]{KKP} where it is verified that their method works perfectly for the twistor 3-fold, predicting both the KKV conjecture and the refinement conjectured in \cite{KKP}.} Instead, one can use their conjectural universal formula
\beq{GV}
\mathop{\sum_{g\ge0}}_{\beta\ne0}
N_{g,\beta}(X)u^{2g-2}v^\beta
\ =\ \mathop{\sum_{g\ge0}}_{\beta\ne0}
n_{g,\beta}(X)\sum_{d>0}\big(2\sin(du/2)\big)^{2g-2\,}\frac{v^{d\beta}}d
\eeq
as the \emph{definition} of the BPS numbers $n_{g,\beta}(X)$. The conjecture is then that $$n_{g,\beta}(X)\in\Z$$ for all $g,\beta,X$, and that for fixed $\beta$ the $n_{g,\beta}(X)$ vanish for $g\gg0$.

We will find it convenient to express Gromov-Witten invariants in terms of BPS numbers in stating the KKV conjecture. It is important to note that \eqref{GV} is an upper triangular linear relationship between the two sets of invariants, with 1s on the diagonal. Hence either set of numbers determines the other; they are equivalent data. \medskip

For instance, for $\beta$ irreducible\footnote{This means there is no decomposition $\beta=\beta_1+\beta_2$ with $\beta_i$ both containing nonzero curves.} we find $n_{0,\beta}(X)=N_{0,\beta}(X)$, but
\beq{AM2}
n\_{0,2\beta}(X)\,=\,N_{0,2\beta}(X)-\frac1{2^3}N_{0,\beta}(X)\,,
\eeq
just as might be expected from the Aspinwall-Morrison formula \eqref{AspMor}. In other words, any genus 0 curve in class $\beta$ contributes $1/2^3$ via double covers
to the genus 0 Gromov-Witten invariant in class $2\beta$, and once we subtract this we expect an integer count of non-multiply covered curves.

\section{K3 surfaces and Noether-Lefschetz loci} \label{K3sec}

For $S$ a projective K3 surface, \eqref{vd} gives $\vd(M_g(S,\beta))=g-1$, even though embedded genus $g$ curves in K3 surfaces move in $g$-dimensional linear systems. The discrepancy is accounted for by the one dimensional
\beq{SR}
H^{0,2}(S)\,\cong\,\C.
\eeq
Deformations of curves are obstructed if we deform $S$ so that the class\footnote{We do not distinguish between $\beta\in H_2(S)$ and its Poincar\'e dual $\beta\in H^2(S)$.} $\beta$ picks up a nonzero component in \eqref{SR} so that it is no longer of Hodge type (1,1). After such deformations, there are no curves in class $\beta$, and the Gromov-Witten invariants of $S$ must vanish by deformation invariance. More directly, \eqref{SR} gives a trivial piece $\O_{\Mb_g(S,\beta)}$ of the obstruction sheaf, and this forces the virtual cycle to vanish.

One can in fact remove this trivial piece \eqref{SR} of the obstruction theory, and get a new 
{\em reduced} virtual cycle
\beq{reddef}
\big[\Mb_g(S,\beta)\big]^{\red}\in A_g\big(\Mb_g(S,\beta)\big)\To H_{2g}\big(\Mb_g(S,\beta)\big)
\eeq
of dimension $g$ whose deformation invariance holds only for $S$ inside the \emph{Noether-Lefschetz locus}\footnote{Here we are deliberately vague about the markings required to make precise sense of this moduli space. For details see \cite{MP}.}
\beq{NL}
N\!L_{\beta}=\big\{S\ | \ \beta\in H^{1,1}(S)\big\}
\eeq
of K3 surfaces for which $\beta$ has type $(1,1)$. \medskip

Alternatively, a three dimensional point of view can be taken. 
The Noether-Lefschetz locus \eqref{NL} is a divisor, given locally by the vanishing of $\int_\beta\sigma_S$, where $\sigma_S$ is a symplectic form on $S$. Therefore we can deform the moduli point
$[S]$ out of it. A generic holomorphic disc $\Delta$ through $[S]\in N\!L_\beta$ will always intersect $N\!L_\beta$.
\begin{figure}[h]
\includegraphics[width=10cm, bb=0 460 600 800]{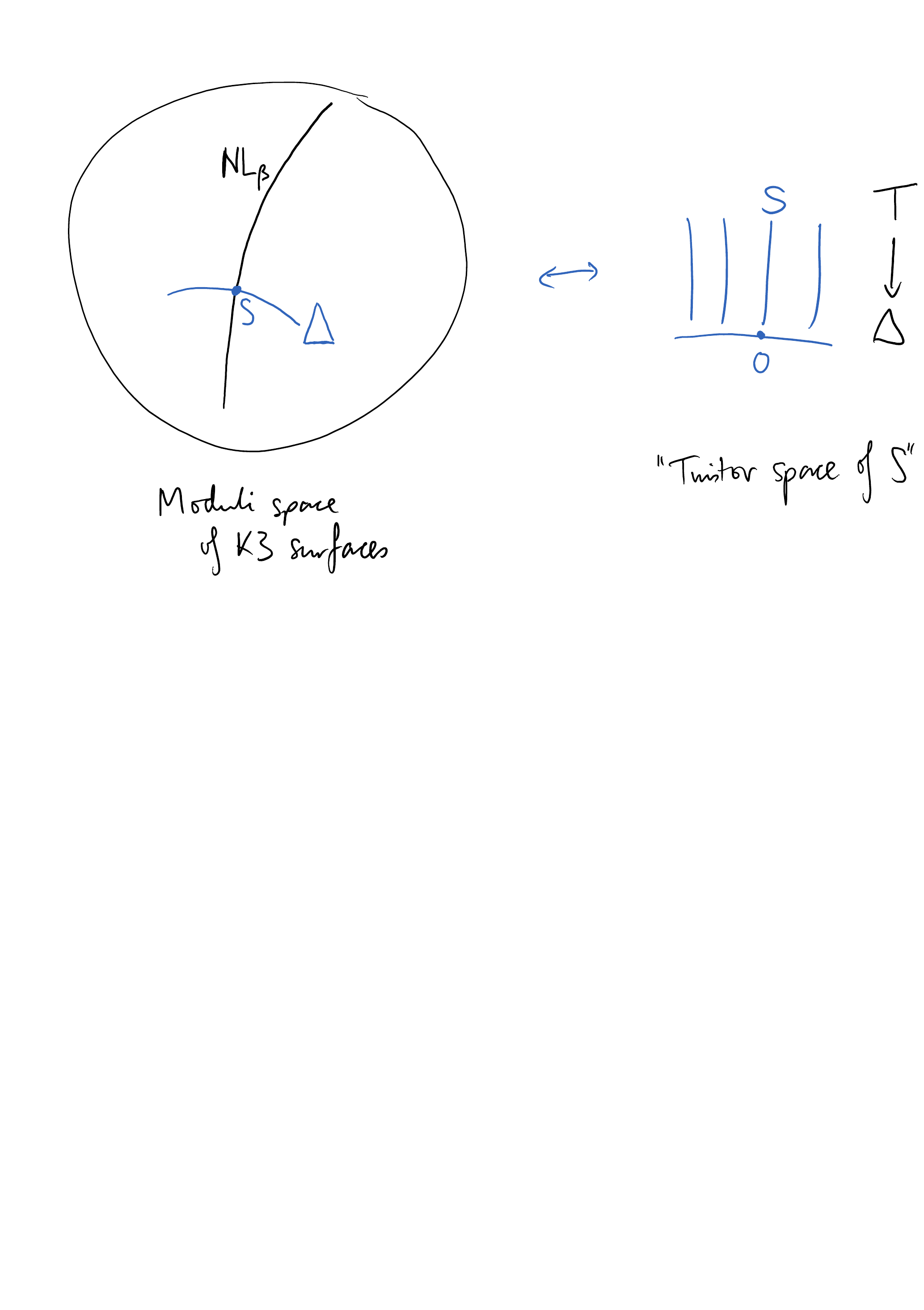}
\end{figure}

Such a disc corresponds to a K3-fibration
\beq{T}
T\To\Delta
\eeq
with central fibre $\iota\colon S\into X$. We take $\Delta$ to be transverse to $N\!L_\beta$ and all $N\!L_\gamma$ for $\gamma$ of smaller degree than $\beta$.  Here $T$ stands informally for {\em twistor space}. In \cite{PTKKV}, algebraic approximations to the real twistor space are used (in order to be able to use the algebraic theory of stable pairs) and $\Delta$ is replaced by a quasi-projective curve. 
Such technicalities add  to the complexity of the paper. For this survey,
 it is safe to think of $T$ as a piece of the actual twistor space of $S$, so that curves in the fibres of $T\to\Delta$ can only appear in the central fibre $S$.

Therefore the moduli spaces
$$
\Mb_g(T,\iota_*\beta)\ =\ \Mb_g(S,\beta)
$$
coincide as sets, and in fact as Deligne-Mumford stacks by a simple deformation theory argument. 
But they have different obstruction theories of different virtual dimensions. 
Since $T$ is a 3-fold whose canonical bundle is trivial on $\iota_*\beta$,
$$\vd(M_g(T,\iota_*\beta))\,=\,0$$
by \eqref{vd}. Thus we can finally define the invariants that this paper is concerned with as
\beq{defn}
N_{g,\beta}(S)=\int_{[M_g(T,\iota_*\beta)]^{\vir}}1.
\eeq
It is not immediately obvious that \eqref{defn} depends
 only on $S$ and not $T$, but by a simple comparison of obstruction theories, one can also express 
\eqref{defn} in terms of the reduced cycle \eqref{reddef} of $S$ as
$$
N_{g,\beta}(S)=\int_{[M_g(S,\beta)]^{\mathrm{red}}}(-1)^g\lambda_g.
$$
Here $\lambda_g=c_g(\mathbb E_g)$, where $\mathbb E_g$ is the Hodge bundle whose fibre over the map $f\colon C\to S$ is $H^0(C,\omega_C)$. \medskip

The upshot is that $N_{g,\beta}(S)$ can be thought of as the contribution of $S$ to the Gromov-Witten theory (in a fibre class) of any K3-fibred 3-fold in which it appears as a fibre.

By deformation invariance and the Torelli theorem for K3 surfaces, the Gromov-Witten invariant $N_{g,\beta}(S)$ depends only on $\beta$ through two integers: $\beta^2$ and the divisibility of $\beta$ in $H_2(S,\Z)$.

\section{The Katz-Klemm-Vafa conjecture} \label{KKVsec}

Since the $N_{g,\beta}(S)$ are really 3-fold invariants \eqref{defn}, it makes sense to rewrite them in BPS form via \eqref{GV}. This gives the equivalent invariants $n_{g,\beta}(S)$. The KKV conjecture is not just that $n_{g,\beta}(S)\in\Z$, but that a further miracle occurs:
\beq{KKV1}
n_{g,\beta}(S)\text{ depends only on }\beta^2,\text{ {\bf not} the divisibility of }\beta.
\eeq
If $\beta^2=2h-2$, we may denote $n_{g,\beta}(S)$ by $n_{g,h}$. The KKV formula
then determines all $n_{g,h}$:
\beq{KKV2}
\sum_{g,h\ge0}(-1)^gn_{g,h\!}\left(\!\sqrt z-\frac1{\sqrt z}\right)^{\!\!2g}\!q^h
\ =\ \prod_{n=1}^\infty\frac1{(1-q^n)^{20}(1-zq^n)^2(1-z^{-1}q^n)^2}\,.
\eeq

The key is the remarkable claim \eqref{KKV1}, first conjectured explicitly in \cite{MP}. It says that whatever $n_{g,\beta}(S)$ counts, it does not see multiple curves or the divisibility of $\beta$. To calculate, we may pass to a completely different homology class of the same square but divisibility 1, count there, and get the same answer! For instance in the simplest example \eqref{AM2}, this says that for primitive classes $\beta,\,\gamma$ with $\gamma^2=(2\beta)^2$ we have $$N_{0,2\beta}(S)-\frac1{2^3}N_{0,\beta}(S)\,=\,N_{0,\gamma}(S).$$

Once \eqref{KKV1} is proved, one may work with primitive classes, for which the KKV formula was proved in \cite{MPT}. However, both are proved together in \cite{PTKKV}. \medskip

The formula \eqref{KKV2} of course also implies that the $n_{g,h}$ are integers.
In particular it gives
$$
n_{g,h}=0 \text{ for } g>h,
$$
which is another remarkable property of the BPS formalism. It says that on K3 surfaces we do not count any maps from higher genus $g>h$ curves to an image curve of arithmetic genus $h$: BPS numbers do not count multiple covers or degenerate contributions. When $g=h$ we find $$n\_{h,h}=(-1)^h(h+1),$$ the signed Euler characteristic of the $\PP^h$ linear system of embedded curves in class $\beta$. For $g<h$ the $n_{g,h}$ count partial normalisations of singular embedded curves of geometric genus in the interval $[g,h]$. 
For small $g,h$ they are \vspace{2mm}
\begin{center}
\begin{tabular}{|c||ccccc|}
        \hline
\textbf{}
$n_{g,h}$&    $h= 0$ & 1  & 2 & 3 & 4 \\
        \hline \hline
$g=0$ & $1$ & $24$ & $324$ & 
$3200$ &$25650$  \\
1      &  & $-2$ & 
$-54$ & $-800$  & $-8550$      \\
2      & & & $3$ & 
$88$ & $1401$       \\
3      & &  & 
 & $-4$  & $-126$       \\
4      &  &  & 
 &   & 5       \\
       \hline
\end{tabular}
\end{center}
\vspace{2mm}
Setting $z\mapsto1$ in \eqref{KKV2} restricts to genus 0 invariants and recovers the Yau-Zaslow formula \cite{YZ}. This was first proved for all classes and all multiple covers in \cite{KMPS} using mirror symmetry.

To prove the KKV conjecture we do not attempt to compute the Gromov-Witten invariants directly. Instead we calculate with the closely related curve-counting theory of stable pairs \cite{PT1}. For a more thorough introduction we again refer to \cite{13.5}.

\section{Stable pairs} \label{SP}

A \emph{stable pair} on a nonsingular projective variety $X$ is a pair $(F,s)$,
\begin{itemize}
\item $F$ is a coherent sheaf with 1-dimensional support,
\item $s\in H^0(F)$,
\end{itemize}
 which is stable:
\begin{itemize}
\item $F$ is \emph{pure:} it has no 0-dimensional subsheaves, and
\item $s$ has 0-dimensional cokernel.
\end{itemize} \medskip

Roughly speaking, a stable pair is the data of a \emph{Cohen-Macaulay curve}\footnote{This can be disconnected, reducible or nonreduced, but has no embedded points.} $$C= \text{supp}(F)$$ plus a \emph{0-dimensional subscheme} $Z\subset C$ (the support of coker$(s)$).
So, for instance,
\begin{enumerate}
\item[(i)] $Z=\emptyset\subset C$ corresponds to the stable pair $(\O_C,1)$.
\item[(ii)] A Cartier divisor $Z\subset C$ corresponds to the stable pair $(\O_C(Z),s\_Z)$.
\item[(iii)] An example with a Weil divisor $Z\subset C$ is given by 
$$C=C_1\cup C_2 \ {\text{and}}\ Z=C_1\cap C_2\, ,$$ corresponding to the stable pair $\big(\O_{C_1}\oplus\O_{C_2},(1,1)\big)$.
\end{enumerate}

There is a projective moduli space
$$
P_n(X,\beta)
$$
of stable pairs with curve class and holomorphic Euler characteristic
$$[F]=\beta,\ \ \chi(F)=n.$$
Thus $(\beta,n)$ is equivalent information to $(c_2(F),c_3(F))$. In fact, $\beta=-c_2(F)$ and $1-n$ is the arithmetic genus of $C$ minus the length of coker$(s)$.\medskip

When $X=S$ is a surface then the arithmetic genus $h$ of $C$ depends only on $\beta$ by adjunction, and the correspondence between stable pairs and curves with 0-dimensional subschemes becomes more precise \cite{PT3}:
\beq{hilb}
P_{1-h+n}(S,\beta)\,=\,\Hilb^n\!\big(\mathcal C/\Hilb_{\beta}(S)\big).
\eeq
Here $\mathcal C\to\Hilb_\beta(S)$ is the universal curve over the Hilbert scheme of curves in class $\beta$, and $\Hilb^n$ denotes the relative Hilbert scheme of points on its fibres. \medskip

By contrast, for 3-folds the genus can jump with the number of points (while keeping their difference $\chi(F)$ constant) as the following picture of a family of stable pairs illustrates.
\begin{center}
\includegraphics[width=10cm, bb=60 640 600 800]{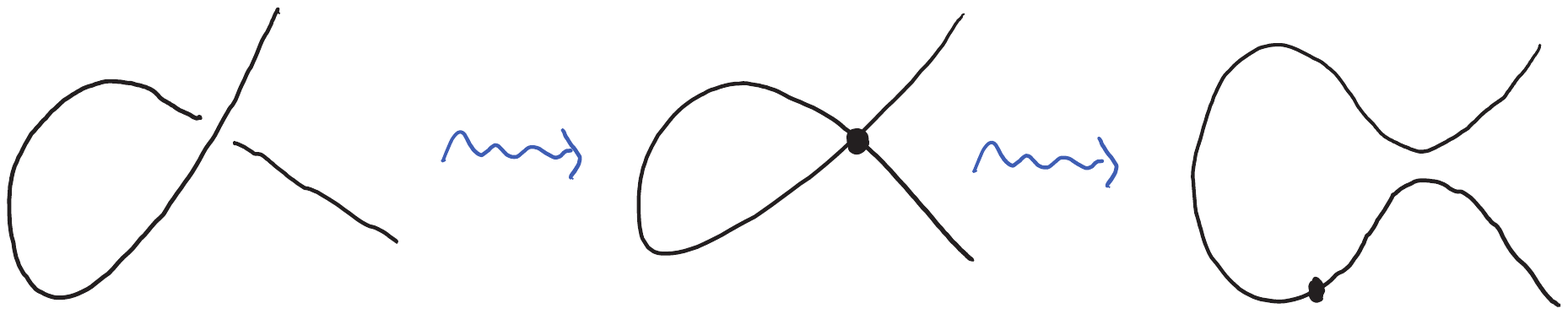}
\end{center}
In terms of the list (i)--(iii) above, this is an example of (i) degenerating to (iii) and then deforming to (ii).

\subsection*{Invariants}
Here it is important that $X$ is a \emph{3-fold}. We consider the complex
$$
I\udot=\{\O_X\Rt{s}F\}
$$
as an object in the derived category of coherent sheaves $D(X)$ of \emph{fixed determinant} $\det I\udot=\O_X$. Its quasi-isomorphism class can be shown to determine the pair $(F,s)$ \cite{PT1}. Its deformations and obstructions are governed by $$\Ext^i(I\udot,I\udot)_0, \quad i=1,2.$$ Their dimensions can jump around, but their difference is the constant
\beq{vdP}
\vd=\int_{\beta}c\_1(X)
\eeq
as before \eqref{vd}. There is a corresponding virtual cycle
$$
[P_n(X,\beta)]^{\vir}\in A_{\vd}\big(P_n(X,\beta)\big)
$$
against which we can integrate to give
\emph{integer} invariants. In our application \eqref{vdP} will be zero, so the invariant is just the degree of the virtual cycle:
\beq{Pndef}
P_{n,\beta}(X)=\int_{[P_n(X,\beta)]^{\vir}}1\ \in\,\Z.
\eeq
This is a different virtual count of curves in class $\beta$, with the genus parameter in Gromov-Witten theory replaced by the closely related $n=\chi(F)$.
The MNOP conjecture \cite{MNOP} states that the two sets of curve-counting invariants -- the rational numbers \eqref{CYGW} defined via stable maps, and the integers \eqref{Pndef} defined via sheaves -- contain the same information.

\section{The MNOP conjecture} \label{MNOPsec}

Let $X$ be a nonsingular projective 3-fold satisfying
the Calabi-Yau condition
$$\int_\beta c\_1(X)=0\, .$$ 
We will discuss the MNOP conjecture as reformulated in  \cite{PT1}
via stable pairs.\footnote{The original MNOP conjecture of \cite{MNOP} relates
 Gromov-Witten invariants to invariants counting ideal sheaves. These latter invariants were conjectured in \cite{PT1} to be equivalent to the stable pair invariants by a certain wall crossing formula. This has now been proved by Bridgeland and Toda \cite{Br, To} in the Calabi-Yau case. The equivalence of ideal sheaf and stable
pairs counting is open for general 3-folds, with some partial results in \cite{StTh}.}

In brief, the MNOP conjecture is the equality
\beq{MNOP}
\exp\Bigg(\sum_{\beta\ne0}\sum_{g\ge0}N_{g,\beta}(X)u^{2g-2}v^\beta\Bigg)\,
=\ 1+\sum_{\beta\ne0}\sum_nP_{n,\beta}(X)q^nv^\beta
\eeq
obtained after substituting $q=-e^{iu}$.

More precisely, the conjecture is first that the generating series of stable pair invariants on the right hand side of \eqref{MNOP} is \emph{the Laurent series of a rational function in $q$ which is invariant\,\footnote{The Laurent series itself need not be invariant under $q\leftrightarrow q^{-1}$. For instance, the series $q-2q^2+3q^3-4q^4+\ldots\ $ does not have this invariance, but its sum is the rational function $$\frac q{(1+q)^2}$$ which \emph{is} invariant under $q\leftrightarrow q^{-1}$.} under $q\leftrightarrow q^{-1}$.} Therefore
the change of variables $q=-e^{iu}$ makes sense, and may be viewed
as the unique analytic continuation from $q=0$ to $q=-1$.

The exponential \eqref{MNOP} turns (the generating series of) connected Gromov-Witten invariants into disconnected invariants: stable pairs is a disconnected theory. Taking logs we can formally define the \emph{connected stable pair invariants} $P^{\;\mathrm{conn}}_{n,\beta}(X)\in \mathbb{Q}$ by the formula
\beq{log}
\sum_{\beta\ne0}\sum_nP^{\;\mathrm{conn}}_{n,\beta}(X)q^nv^\beta=\log\Bigg(
1+\sum_{\beta\ne0}\sum_nP_{n,\beta}(X)q^nv^\beta\Bigg).
\eeq
The MNOP conjecture \eqref{MNOP} can now be written for fixed $\beta$:
\beq{MNOP2}
\sum_{g\ge0}N_{g,\beta}(X)u^{2g-2}\ =\ 
\sum_n P^{\;\mathrm{conn}}_{n,\beta}(X)q^n, \qquad q=-e^{iu}.
\eeq

Because of the  change of variables $q=-e^{iu}$,
it is not sensible to try to interpret the MNOP conjecture at the level of coefficients.\footnote{When put into BPS form the formula is more comprehensible, especially in the irreducible case, see \cite[Section $4\frac12$]{13.5} for a discussion.} The moral is simply that both theories contain the same information, and that there are integers\footnote{In fact the integrality predicted by the MNOP conjecture is the same as that predicted by the Gopakumar-Vafa conjecture, i.e. the $P_{n,\beta}(X)$ are integers if and only if the $n_{g,\beta}(X)$ are integers \cite[Section 3.5]{PT1}.} $P_{n,\beta}(X)$ underlying that rational Gromov-Witten invariants $N_{g,\beta}(X)$.

The most powerful tool we have in understanding Gromov-Witten invariants of threefolds is the following. We refer to \cite{PaPix} for the precise statement.

\begin{thm}[Pandharipande-Pixton \cite{PaPix}] \label{PP}
The MNOP conjecture is true for projective Calabi-Yau 3-folds which can be degenerated to unions of toric varieties.
\end{thm}

\section{Local MNOP for the twistor space} \label{localMNOP}

We can not immediately apply Theorem \ref{PP} to our situation since the threefold $T$ of \eqref{T} is not projective. We have to work hard to prove a \emph{local} form of MNOP for $T$ in the next 5 Sections.

The idea is that in a K3-fibration $X\to C$, curves in fibre classes only appear in the fibres which lie in the Noether-Lefschetz locus. Assuming for simplicity that $C$ is transverse to the relevant Noether-Lefschetz loci, a neighbourhood of each such fibre is well modelled by $T$ and should contribute the invariants of $T$. Applying Theorem \ref{PP} to a projective K3-fibration, then, we would like to show that a combination of the Gromov-Witten invariants of $T$ equals the same combination of \emph{connected} stable pair invariants of $T$. With a judicious choice of K3-fibration we will find the resulting expressions of one set of invariants in terms of the other can be inverted to prove their equality. This will be the local MNOP correspondence we need.

The first part of the above sketch is the {\em Gromov-Witten/Noether-Lefschetz correspondence} of \cite{MP}.

\begin{thm}[Maulik-Pandharipande \cite{MP}] \label{GWNL}
For $X\to C$ a \emph{projective} K3-fibration,
$$
N_{g,\beta}(X)=\sum_{h,m}N_{g,m,h}(S)\cdot N\!L_{m,h,\beta}(X/C)\, .
$$
\end{thm}

\noindent Here the three invariants related by the GW/NL
correspondence are: \begin{itemize}
\item \emph{$N_{g,\beta}(X)$ is the full Gromov-Witten invariant of $X$ in a fibre class $\beta$.}
\item \emph{$N_{g,m,h}(S)$ is the invariant $N_{g,\gamma}(S)$ \eqref{defn} of any K3 surface $S$ and class $\gamma\in H^{1,1}(S,\Z)$ of divisibility $m$ and square $\gamma^2=2h-2$.}
\item \emph{$N\!L_{m,h,\beta}(X/C)$ counts the K3 fibres of $X\to C$ for which a curve class (of square $2h-2$, divisibility $m$ and pushforward $\beta\in H_2(X,\Z))$ becomes of Hodge type $(1,1)$.}
\end{itemize}
More precisely $N\!L_{m,h,\beta}(X/C)$ is an intersection number of $C$ with a Noether-Lefschetz divisor, see \cite{MP} for full details. The result also holds when $C$ does not intersect the Noether-Lefschetz divisors transversely.

To combine this with Theorem \ref{PP} we need an analogous 
{\em Pairs/Noether-Lefschetz correspondence}. This is more complicated since stable pairs count \emph{disconnected} curves, and also since stable pairs in $T$ do not need to lie scheme-theoretically in the central fibre $S$: the multiple covers of Gromov-Witten theory get replaced by scheme-theoretic thickenings of curves in the normal direction to $S\subset T$.

In fact, it is not even obvious how to define a stable pair invariant of $S$ analogous to \eqref{defn}: taking stable pairs on $T$ may give an invariant which depends upon $T$. To get around this, we will deform $T$ to the normal cone of $S\subset T$.

\section{Degeneration} \label{JLi}

Suppose a nonsingular variety $X_t$ degenerates to a variety
$$
X_0=X_1\cup_D X_2
$$
which is a normal crossings divisor in the total space.
Na\"ively, we might expect to find something like
\beq{split}
\big\{\mathrm{curves\ on\ }X\big\}\ \sim\ \big\{\mathrm{curves\ on\ }X_1\big\}
\times\_{\Hilb^d\!D}\big\{\mathrm{curves\ on\ }X_2\big\},
\eeq
where $d$ is the intersection number $\beta\cdot D$ and the map from $\big\{\mathrm{curves\ on\ }X_i\big\}$ to $\Hilb^d\!D$ is a boundary map, intersecting the curve with the divisor $D$.

For stable maps, such a theory originated in \cite{Ruan}
and was developed in algebraic geometry in \cite{LiRelative, junli}. 
We will concentrate on the parallel story for stable pairs \cite{liwu,PT1}. 
To make  \eqref{split} work, we have to avoid the situation where
a component of the curve  falls into $D$. If a component comes close
to falling into $D$, we bubble the target:{\footnote{We have pictured the case relevant to us, where $D$ has trivial normal bundle 
$$N_D\cong\O_D\, .$$ In general $D\times\PP^1$ is replaced by the projective completion $\PP(N_D\oplus\O_D)$ of $N_D$.}}
\begin{center}
\includegraphics[width=12cm, bb=40 610 600 810]{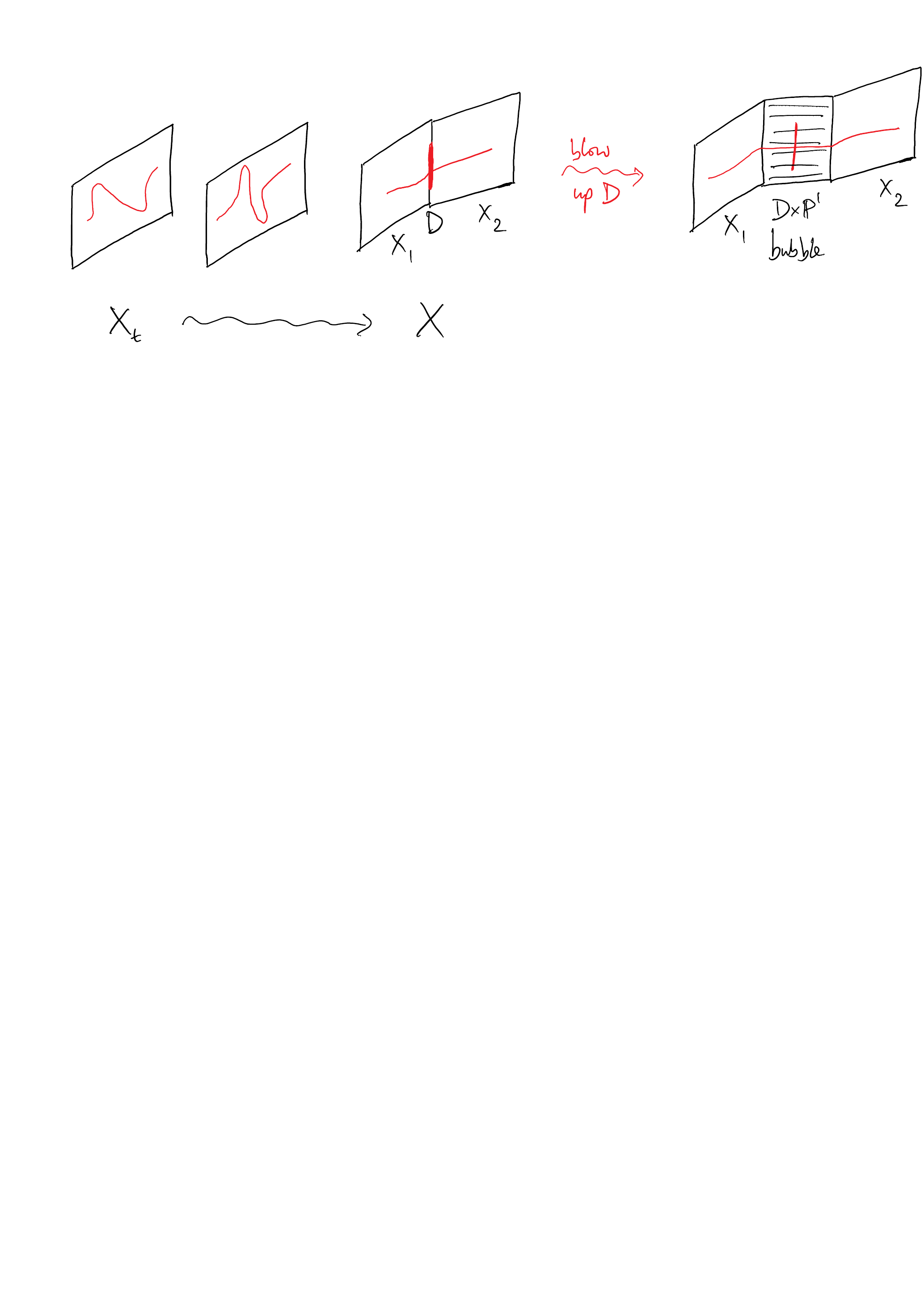}
\end{center}
 To avoid curves falling into the new copy of $D$, we may have to bubble again and again.

In \cite{liwu} J. Li and B. Wu construct \emph{compact} moduli spaces\footnote{These moduli spaces of \emph{relative} stable pairs are what play the role of $\{$curves on $X_i\}$ in \eqref{split}.} $P_{n_i}(X_i/D,\beta_i)$ of stable pairs on $X_i$ \emph{relative to $D$} which admit boundary maps 
$$
P_{n_i}(X_i/D,\beta_i)\To\Hilb^{d\!}D
$$
and virtual cycles $[P_{n_i}(X_i/D,\beta_i)]^{\vir}$ such that
$$
[P_n(X_0,\beta)]^{\vir}\ =\!\!
\mathop{\sum_{\beta_1+\beta_2=\beta}}_{n_1+n_2=n+d}[P_{n\_1}(X_1/D,\beta_1)]^{\vir}
\times_{\Hilb^{d\!}D}[P_{n\_2}(X_2/D,\beta_2)]^{\vir}
$$
is a specialisation of the cycles $[P_n(X_t,\beta)]^{\vir}$. We now explain these statements in a more detail.
\begin{center}
\includegraphics[width=13cm, bb=40 590 600 780]{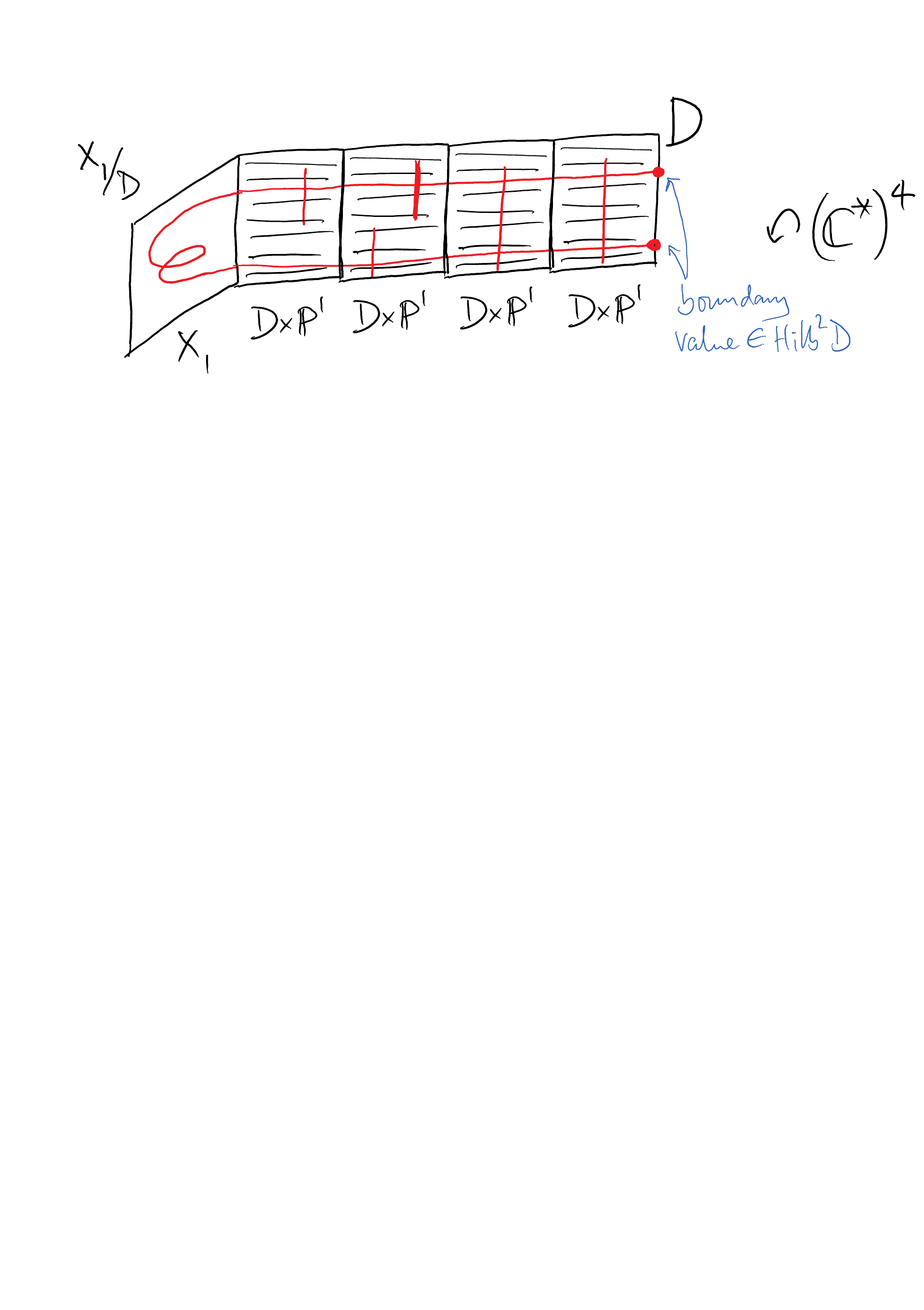}
\end{center}

The moduli space  $P_{n}(X_i/D,\beta_i)$ parameterises isomorphism classes of stable relative pairs
\begin{equation}\label{vyq}
\O_{X_i[k]} \stackrel{s}{\rightarrow} F
\end{equation}
on some $X_i[k]$, the $k$-step degeneration of $X_i$
along $D$. Here, $X_i[k]$ is $X_i$
with $k\ge0$ bubbles attached along $D$ \cite{LiRelative}, and
 $F$ is a sheaf on $X_i[k]$ with $\chi(F)=n_i$ whose support pushes down to $\beta_i\in H_2(X_i,\Z).$
The stability conditions for the data are more complicated
in the relative geometry, but ensure that the stable pair joins correctly across the creases of $X_i[k]$ and intersects the last copy of $D$ transversally:
\begin{enumerate}
\item[(i)] $F$ is pure with finite locally free resolution,
\item[(ii)] the higher derived functors of the
restriction of $F$ to the singular loci of $X_i[k]$, and to the final copy of $D$, vanish,
\item[(iii)] the section $s$ has 0-dimensional cokernel supported
away from the singular loci of $X_i[k]$,
\item[(iv)] the pair \eqref{vyq} has only finitely many automorphisms covering
the automorphisms of $X_i[k]/X_i$.
\end{enumerate}
Relative stable pairs are isomorphic if they differ by an element of
$$
\Aut(X_i[k]/X_i)\,=\,(\C^*)^k.
$$
Condition (iv) ensures that we do not insert unwanted bubbles and that
$P_{n_i}(X_i/D,\beta_i)$ is a Deligne-Mumford stack. By conditions (ii) and 
(iii), the moduli space admits a boundary map to $\Hilb^{d\!}D$ by restriction.

\section{Deformation to the normal cone} \label{ncone}

We apply the relative theory of Section \ref{JLi} to calculate 
the stable pair invariants of the {twistor 3-fold} $T$ by deforming $T$ to the normal cone of the central fibre $S$.
\begin{center}
\includegraphics[width=8cm, bb=50 600 450 810]{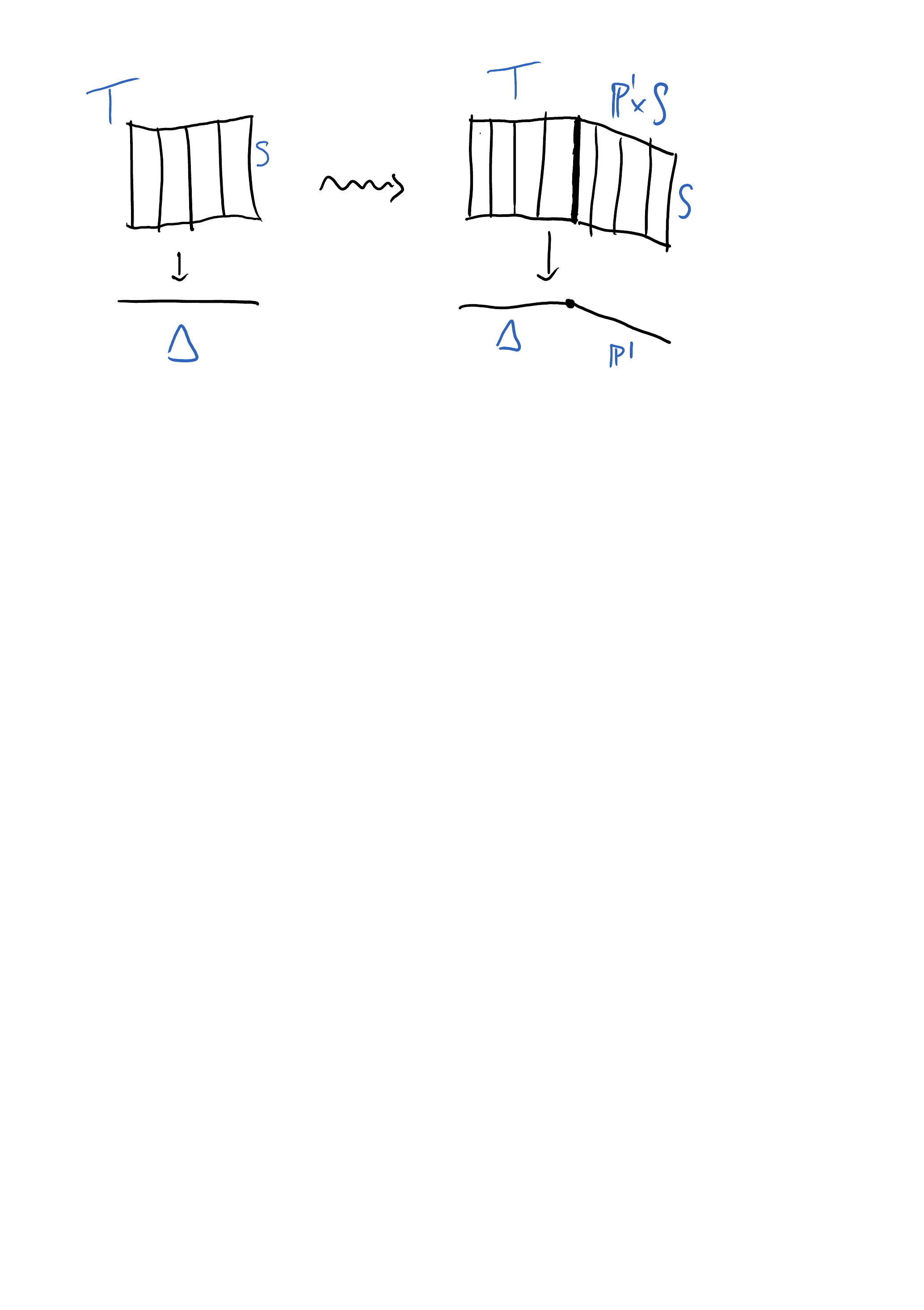}
\end{center}
We may compute the stable pair invariants of $T$ via the 
degeneration formula in terms of  
the relative theories of $T$ and $S\times\PP^1$  attached along $S$.

Since we are in a fibre class $\beta$, matters simplify. 
Curves are not allowed in the creases in the relative theory, 
so they cannot lie in the central fibre $S$ of $T$ (along which $S\times\PP^1$ is attached). 
By construction, there are no curves in class $\beta$
 in the other fibres of $T$. Hence, {\em all} the relevant curves  lie in 
$S\times \PP^1$ and its bubbles.\footnote{In the notation of Section \ref{JLi}, we have $d=\beta\cdot S=0$.}
\begin{center}
\includegraphics[width=12cm, bb=200 700 600 810]{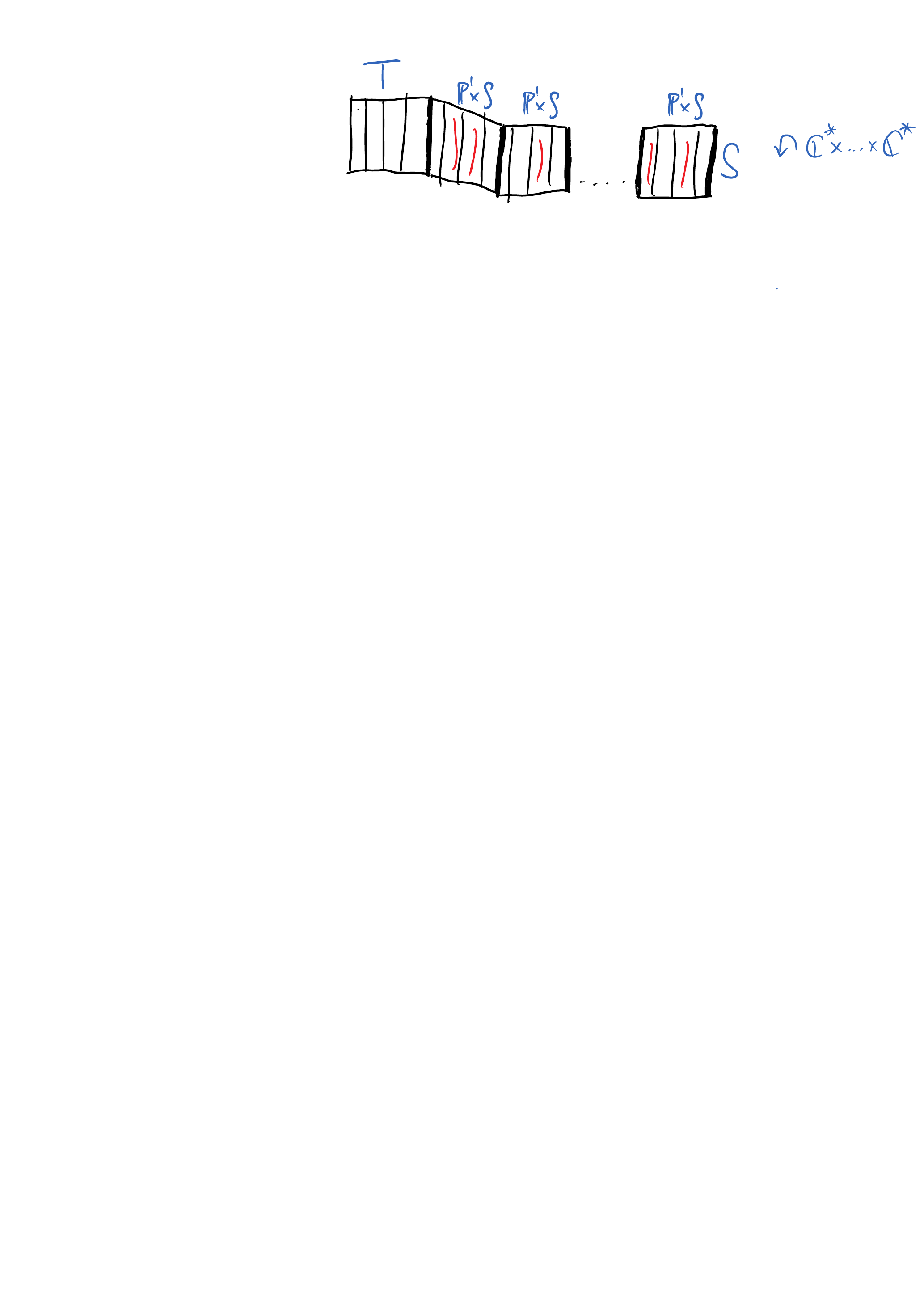}
\end{center}
We have made progress: $S\subset T$ 
has been replaced by the trivial fibration $S\times\PP^1$ (together with its bubbles). We see
\begin{itemize}
\item the answer depends only on $S$ after all (this involves a careful study of the perfect obstruction theory, which \emph{a priori} depends on $T$),
\item we can use $\C^*$-localisation.
\end{itemize}
The advantages come at the expense of having passed to the relative theory with bubbles and $(\C^*)^k$-automorphisms, but these turn out to be manageable.

\section{Connected stable pairs theory} \label{conn}

We have degenerated the moduli space of stable pairs on $T$ to the moduli space
\beq{rub}
P_n\left( S\times R,\beta\right)
\eeq
of stable pairs on the {\em rubber} geometry 
$$S \times R =  S\times\PP^1 \, \big/\, (S\times\{0\}) \cup (S\times\{\infty\})\, $$
modulo identification of pairs which differ by the action of $\C^*$. It is a compact space.
\begin{center}
\includegraphics[width=12cm, bb=40 680 550 820]{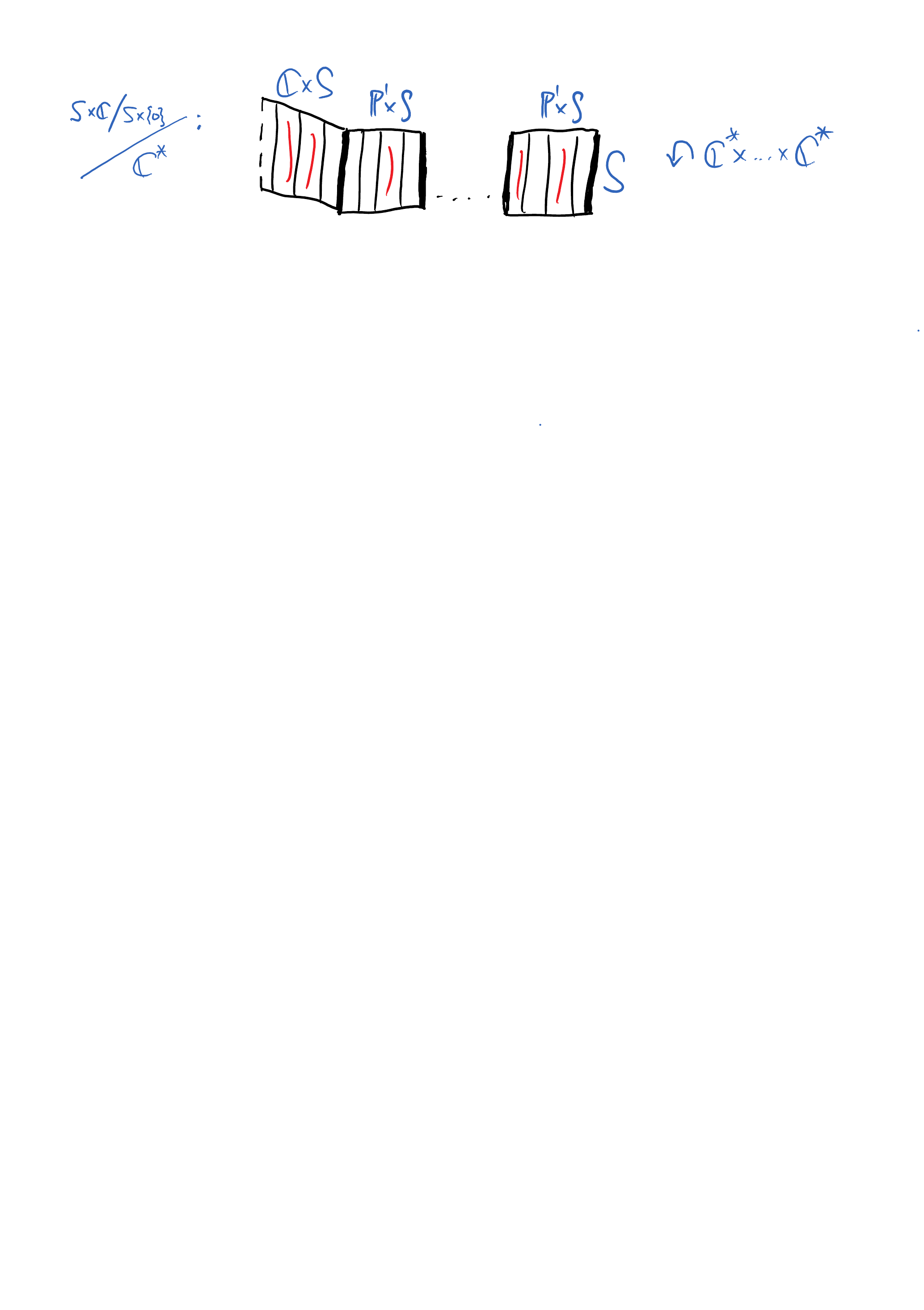}
\end{center}
Na\"ively \eqref{rub} does not see $T$, but in fact its perfect obstruction theory -- induced by the relative theory applied to the degeneration to the normal cone -- does. In particular, it differs from the standard perfect obstruction theory on
\eqref{rub}, so the virtual cycles differ too.

Analysing their difference and using an induction on the number of bubbles gives the following surprising formula for the invariants
$$
1+\sum_{n,\beta\ne0}P_{n,\beta}(T)q^nv^\beta\ =\,\exp\Bigg(\sum_{n,\beta\ne0}P^{\red}_{n,\beta}
\big( S\times R)q^nv^\beta\Bigg).
$$
Here, $P^{\red}_{n,\beta}\big(S\times R \big)$ is the \emph{reduced} stable pair invariant of the rubber geometry $S\times R$, 
defined by removing the trivial $H^{0,2}(S)=\C$ piece from the natural perfect obstruction theory just as in Section \ref{K3sec}.

Put differently, we have found that
$$
P^{\mathrm{\;conn}}_{n,\beta}(T)\,=\,P^{\red}_{n,\beta}\big( S\times R \big),
$$
which is clearly not dependent on $T$. By deformation invariance of these invariants as $S$ moves within $N\!L_\beta$, they depend only on $\beta^2=2h-2$ and the divisibility $m$ of $\beta$. Therefore, we may write them both as $$R_{n,m,h}(S).$$It is for this connected theory that we can prove a Pairs/Noether-Lefschetz correspondence.

\section{Pairs/Noether-Lefschetz correspondence} \label{PNL}

If $X\to C$ is a projective K3-fibration for which $C$ intersects the relevant Noether-Lefschetz divisors transversally, the same argument as above (deformation to the normal cone of the finitely many fibres in the Noether-Lefschetz locus) gives the following Pairs/Noether-Lefschetz correspondence. 

\begin{thm}[\cite{PTKKV}] \label{SPNL} If the generic Picard rank of
the K3 fibers of $X$ is at least 3,
$$
P^{\mathrm{\;conn}}_{n,\beta}(X)\ =\ \sum_{h,m}R_{n,m,h}(S)\cdot N\!L_{m,h,\beta}(X/C).
$$
\end{thm}
Here, as in Theorem \ref{GWNL}, $N\!L_{m,h,\beta}(X/C)$ is an intersection number of $C$ with a Noether-Lefschetz divisor. It counts the fibres of $X\to C$ for which a curve class of square $2h-2$, divisibility $m$ and push-forward $\beta\in H_2(X,\Z)$ becomes of Hodge type (1,1).

Theorem holds more generally without the transversality assumption. 
If there are only finitely many fibres in the Noether-Lefschetz locus then the fact that all contributions to the above theorem are local to these fibres means that we may perturb $C$ locally analytically to intersect the Noether-Lefschetz divisors transversally. If $C$ lies entirely within a Noether-Lefschetz divisor then we add a very positive curve $C'$ in the moduli space of K3 surfaces and then deform to a smooth curve $C''$ which intersects the Noether-Lefschetz divisors transversally. The degeneration formula for the K3-fibration over $C''$ degenerating to the K3-fibration over $C\cup C'$ then gives the result.

By the technical assumption of generic Picard rank at least 3 for the K3
fibers of $X$,  the perturbations in the
degenerate cases above are easily realizable. Of course the
result is expected to hold without the Picard rank assumption, just
as in the Gromov-Witten case.

\section{Invertibility} \label{invert}

We can now apply the GW/NL correspondence of Theorem \ref{GWNL}, 
the Pairs/Noether-Lefschetz correspondence of Theorem \ref{SPNL},
 and the MNOP correspondence of Theorem \ref{PP} 
to a nonsingular anticanonical divisor
\begin{equation}\label{jzjz}
X\ \subset\ \Bl_p(\PP^2\times\PP^1)\times\PP^1,
\end{equation}
considered as a K3-fibration via projection to the final $\PP^1$ factor.

The result is that the MNOP formula \eqref{MNOP2} for $X$,
$$
\sum_{g\ge0}N_{g,\beta}(X)u^{2g-2}\ 
=\ \sum_nP^{\;\mathrm{conn}}_{n,\beta}(X)q^n, \qquad q=-e^{iu},
$$
is expressed as the equality of a linear combination of terms
\beq{GWT}
\sum_g N_{g,m,h}(S)u^{2g-2},
\eeq
for different $m$ and $h$, and the \emph{same} linear combination of terms
\beq{SPT}
\sum_n R_{n,m,h}(S)q^n.
\eeq
The particular 3-fold $X$ in \eqref{jzjz}
was chosen to make these linear relations invertible. 
Hence, it follows that  \eqref{GWT} and \eqref{SPT} are equal for all $m,h$,
\beq{MNOP-ST}
\sum_g N_{g,m,h}(S)u^{2g-2}\ =\ \sum_n R_{n,m,h}(S)q^n, \qquad q=-e^{iu}.
\eeq
This gives the local MNOP correspondence for $S$ that we sought. That is, the MNOP conjecture is true for the local threefold $T$.

\section{Pairs on $S\times\C$} \label{EC}

We have turned the KKV conjecture into a question
 entirely about (reduced) stable pairs on the rubber geometry 
$S \times R$. By rigidification techniques, further degeneration, 
and finally $\C^*$-localisation, we end up needing only to calculate the $\C^*$-localised \emph{reduced} stable pair invariants of $S\times\C$. \medskip

The latter involves calculating integrals over the moduli spaces
$$
P_n(S\times\C,\beta)^{\C^*}
$$
of $\C^*$-fixed stable pairs. Such pairs are supported set-theoretically on 
$$S\times\{0\} \subset S\times \C\, ,$$
 but can be thickened infinitesimally out of $S\times\{0\}$ into $S\times\C$. These thickenings replace the multiple covers of Gromov-Witten theory, which all lie in $S$.

However, one advantage stable pair theory has over Gromov-Witten theory
is its \emph{symmetric obstruction theory} on Calabi-Yau 3-folds.\footnote{In fact, the moduli space $P_n(S\times\C,\beta)$ is \emph{locally the critical locus of a holomorphic function on a nonsingular space}.} As a consequence,
\beq{dual}
\textit{vector fields on $P_n(S\times\C,\beta)$ are dual to obstructions.} \eeq
We have already used this fact once. The vector field $\partial_x$ on $S\times\C_x$ that translates in the $\C_x$-direction is dual to the trivial piece of the obstruction theory \eqref{SR} that we remove to get the \emph{reduced} obstruction theory of $$P_n(S\times\C_x,\beta)\, .$$ We will now use it a second time, on a related vector field.

\subsection*{Second vector field}
Consider the flow outward from
\beq{incl}
P_n(S\times\C_x,\beta)^{\C^*}\subset\,P_n(S\times\C_x,\beta)
\eeq given by \emph{pulling apart the {\em last layer} of thickening of a stable pair}:\medskip
\begin{center}
\includegraphics[width=10cm, bb=40 500 600 740]{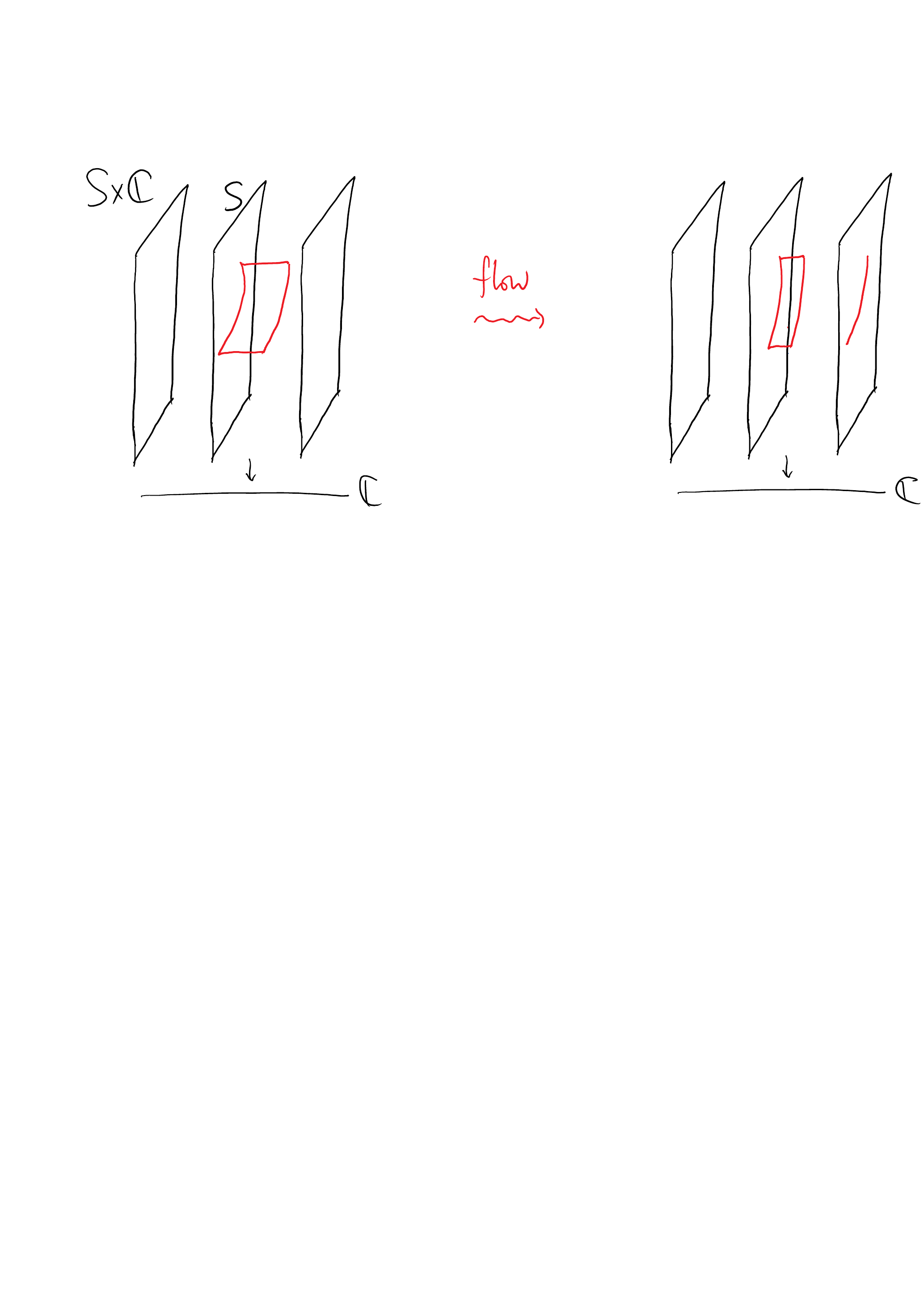}
\end{center}
More precisely, work on a component $P[d]$ of $P_n(S\times\C_x,\beta)^{\C^*}$ in which all stable pairs are supported on 
$$dS=S\times\Spec\C[x]/(x^d)$$ but not on $(d-1)S$. Then, there is a map
\beq{flow}
P[d]\times\C_t\To P_n(S\times\C_x,\beta)
\eeq
which at $t=0$ restricts to the inclusion \eqref{incl}, and for $t\ne0$ takes each stable pair to one supported on $(d-1)S\sqcup S_t$. The precise details are given in \cite[Section 5]{PTKKV}, but the local model in the $\C_x$-direction is the family $\O_Z$, where
\beq{Z}
Z=\big\{x^{d-1}(t-x)=0\big\}\subset\C_x\times\C_t.
\eeq
This is a flat family of schemes over $\C_t$ with central fibre the $d$-times thickened point $\{x^d=0\}\subset\C_x$ at $t=0$, and general fibre
$$\{x^{d-1}=0\}\sqcup\{x=t\}$$ over $t\ne0$.
 In particular, the {\em centre of mass} of the subscheme moves to the fibre
\beq{com}
S_{\frac td}\subset S\times\C
\eeq
at time $t$. \medskip

We only need the flow to first order in $t$. Differentiating \eqref{flow} at $t=0$ gives a vector field
\beq{vfield}
v\in\Gamma\Big(T_{P_n(S\times\C,\beta)}\big|_{P[d]}\Big).
\eeq
That is, it is a $P_n(S\times\C,\beta)$-vector field, but only on the subscheme $P[d]$. By \eqref{com} it moves the centre of mass of the subscheme by the vector field
\beq{comv}
\partial_t/d. 
\eeq

\section{Vanishing result} \label{zero}

It turns out that, to first order in $t$, the above deformation only sees how the centre of mass deforms. That is -- surprisingly, perhaps -- on basechange to $\Spec\C[t]/(t^2)\subset\C_t$, the deformation \eqref{Z} becomes the same as the first order deformation given by moving the whole scheme $Z$ along the vector field $\partial_t/d$ \eqref{comv}. In particular,
\beq{texty}
\textit{the vector field $v$ equals $\frac1d\partial_t$ on $d$-times uniformly-thickened stable pairs.}
\eeq
Here, we say a stable pair is \emph{uniformly-thickened} if it takes the form
\beq{thick}
(F,s)\ =\ (F_0,s_0)\otimes\_\C\frac{\C[x]}{(x^d)}
\eeq
for some stable pair $(F_0,s_0)$ on $S$.

\includegraphics[width=12cm, bb=20 540 590 820]{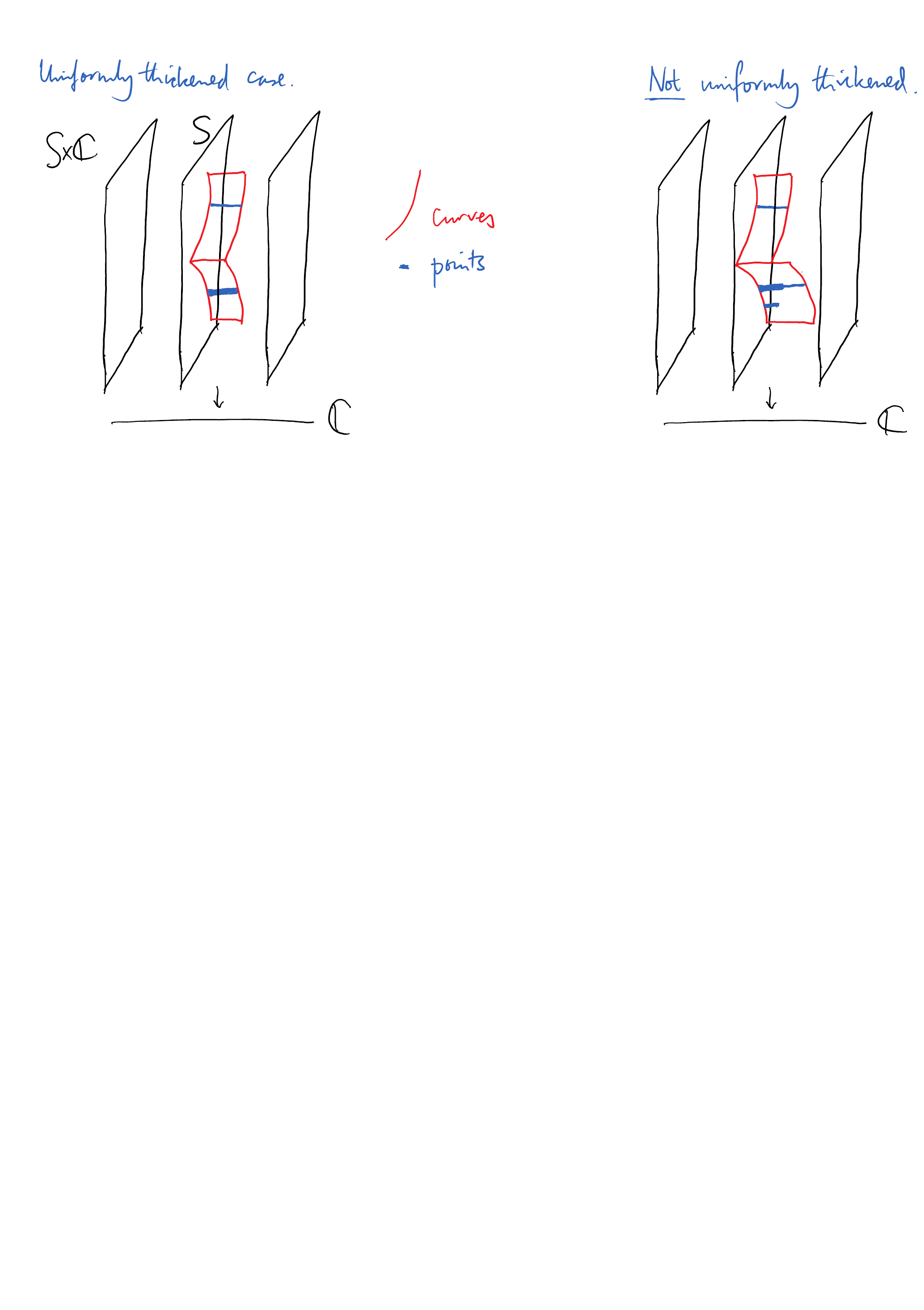}

By \eqref{texty} $v$ is linearly dependent on $\partial_t$
along a component of $$P_n(S\times\C,\beta)^{\C^*}$$ if and only if all pairs in the component are uniformly-thickened in the $\C$ direction.

For $\C^*$-fixed stable pairs which are \emph{not} uniformly-thickened in this way, the vector field $v$ acts as different vector fields $\partial_t/k$ at different points of their support, so it is not globally proportional to any one $\partial_t/k$.

By the duality \eqref{dual}, this means that for non-uniformly-thickened stable pairs, there is \emph{another} trivial piece of the obstruction theory. That is we get a surjection from the obstruction sheaf to $\O_{P[d]}$ dual to $v$ that is not proportional to the original one \eqref{SR}. This forces the reduced virtual cycle to be zero. In particular we get a vanishing result.
 
\begin{prop} \label{vanish}
The only components of $P_n(S\times\C,\beta)^{\C^*}$ which contribute are those consisting of \emph{uniformly-thickened} stable pairs.
\end{prop}

\section{Multiple covers} \label{mult}

Proposition \ref{vanish} means we can really calculate the full $\C^*$-localised stable pairs theory of $S\times\C$ (and so, by our previous results, of $T$ and the K3-fibration $X$).

The moduli space of $d$-times uniformly-thickened stable pairs on $S\times\C$ is just
\beq{/d}
P_{n/d}(S,\beta/d).
\eeq
Given a pair $(F_0,s_0)$ in the above moduli space the corresponding thickened pair in $P_n(S,\beta)^{\C^*}$ is given by \eqref{thick}.

This isomorphism does not preserve the relevant perfect obstruction theories, but the two can be compared. Therefore we have reduced everything to a calculation over a moduli space \eqref{/d} of stable pairs supported \emph{scheme theoretically on $S$}.

Such integrals were considered in \cite{KT1,KT2} in connection with the G\"ottsche conjecture and its virtual extensions. We can use
\eqref{hilb},
\beq{hilb2}
P_{1-h+n}(S,\gamma)\,=\,\Hilb^n\!\big(\mathcal C/\Hilb_{\gamma}(S)\big),
\eeq
which is shown in \cite{KT1} to be an isomorphism not just of schemes, but of \emph{schemes with perfect obstruction theory.} Here we take the reduced obstruction theory on the left, and on the right a natural one given by equations. Namely, the relative Hilbert scheme $\Hilb^n\!\big(\mathcal C/\Hilb_{\gamma}(S)\big)$ lies in the smooth space\footnote{Here, matters are simpler because $S$ is simply connected. For more general surfaces, $\Hilb_\gamma(S)$ need not be smooth, and has to be embedded in the smooth space $\Hilb_{\gamma+A}(S)$, for an ample divisor $A\gg0$.}
\beq{amb}
\PP(H^0(L))\times\Hilb^n(S)
\eeq
of curves in $S$ (in class $c_1(L)=\gamma$) and points in $S$. It is the incidence variety
$$
\big\{(s,Z)\colon s|_Z=0\big\}
$$
of points $Z$ which lie on the curve $\{s=0\}$.
As such it is cut out by the section $s|_Z$ of the tautological vector bundle $E$ whose fibre over $(s,Z)$ is $H^0(L|_Z)$. This description endows \eqref{hilb2} with a perfect obstruction theory whose virtual cycle is
$$
c_{\mathrm{top}}(E)\,=\,[P_{n+1-h}(S,\gamma)]^{\red}
$$
when pushed forward to \eqref{amb}. This allows one to compute the reduced virtual cycle in terms of tautological integrals to which one can apply \cite{EGL}. The result is universal formulae in the topological numbers $\beta^2$ and $n$, but
$$
\textit{not the divisibility of $\beta$.}
$$
Translated through the MNOP correspondence, 
this is what gives the required independence of BPS numbers on the divisibility $m$ \eqref{KKV1}. 

In terms of stable pairs only, the resulting multiple cover formula is the following. Note its simplicity compared to multiple cover formulae in higher genus Gromov-Witten theory.

\begin{prop}\label{fzzq}
The generating series
$$
P^{\red}_{d\beta}(q)=\sum_{n\in\Z}P^{\red}_{n,d\beta}(S\times\C)q^n
$$
of reduced stable pair invariants of $S\times\C$ is the Laurent series of a rational function of $q$. It satisfies the multiple cover formula
$$
P^{\red}_{d\beta}(q)=\sum_{k|d}\frac1kP^{\red}_{\gamma}\big(-(-q)^k\big),
$$
for any primitive class $\gamma$ with the same square as $d\beta/k$.
\end{prop}

Together with \eqref{MNOP-ST}, 
Proposition \eqref{fzzq} gives the insensitivity of BPS numbers to divisibility of the curve class \eqref{KKV1} and allows us to compute with only primitive $\beta$. 

We may further deform to the case where $\beta$ is irreducible. Then, 
the moduli space of stable pairs is nonsingular and the invariants become signed Euler characteristics. Using sheaf theory on K3 surfaces -- a heavily developed and well-behaved subject due to the tight constraints of holomorphic symplectic geometry -- the relevant calculations here were done by Kawai and Yoshioka \cite{KY}, see also \cite{MPT,PT3}. This completes the proof of the KKV formula \eqref{KKV2}.

\bibliographystyle{halphanum}
\bibliography{references}

\bigskip \noindent {\tt{rahul@math.ethz.ch}} \medskip

\noindent Departement Mathematik \\
\noindent ETH Z\"urich \\
\noindent 8092 Z\"urich \\
\noindent Switzerland

\bigskip \noindent {\tt{richard.thomas@imperial.ac.uk}} \medskip

\noindent Department of Mathematics \\
\noindent Imperial College London\\
\noindent London SW7 2AZ \\
\noindent United Kingdom

\end{document}